
\input amstex.tex

 \input amsppt.sty \nologo \NoBlackBoxes
 \magnification=\magstep1

    \topmatter

 \title Picard Group of the Moduli Spaces of $G$--Bundles
 \endtitle
 \author Shrawan Kumar and M. S. Narasimhan \endauthor
 \address S.K.: Department of Mathematics, University of North Carolina,
Chapel Hill, NC 27599--3250, USA\endaddress
\address M.S.N.: ICTP, P. O. Box 586, Strada Costiera 11, Miramare,
I-34100
Trieste, ITALY  \endaddress


    \endtopmatter
  \document
\vskip6ex
\head Introduction \endhead

 Let $G$ be a simple simply-connected connected complex affine
 algebraic group and let $C$ be a smooth irreducible
 projective curve  of genus $\geq 2$ over the field of complex numbers $\Bbb
 C$.
 Let $\frak M$ be the moduli space of semistable
 principal $G$-bundles on $C$ and let Pic $\frak M$
be its Picard group, i.e., the group of isomorphism classes of algebraic line
bundles on $\frak M$. Following is our main result (which generalizes a result
of Drezet-Narasimhan for $G=$ SL$(N)$ \cite{DN} to any $G$).

 \proclaim{(A) Theorem} With the notation as above, {\rm Pic}$\,(\frak
M)\approx \Bbb Z$.
 \endproclaim

A more precise result is obtained in Theorem (2.4) together with Theorem (4.9)
(see also Question 4.13).

 We use the above result and a result of Grauert-Riemenschneider  to prove the
following second main result of this
 paper.

 \proclaim{(B) Theorem} The dualizing sheaf $\omega$ of the
 moduli space $\frak M$ is locally free.
In particular, $\frak M$ is a Gorenstein variety.

Further,  for any finite dimensional representation $V$ of $G$ ,
 $H^i(\frak M,\Theta(V))=0$, for
 all ~ $i>0$, where  $\Theta(V)$ is the theta bundle on the moduli space $\frak
M$. In particular,
 $$
 \Cal X(\frak M,\Theta(V)) = \dim H^0(\frak M,\Theta(V)),
 $$
 where $\Cal X$ is the Euler-Poincare characteristic.
\endproclaim

In fact we have a sharper  result than the above (cf. Theorem 2.8).
\vskip2ex

We make essential use of the flag variety $X$ associated to the affine
Kac-Moody group corresponding to $G$, which parametrizes an algebraic family of
$G$-bundles on $C$, and the fact that Pic $X \simeq {\Bbb Z}$. We also need to
make use of the explicit construction of the moduli space $\frak M$ via GIT.
\vskip 3ex
\noindent
{\bf Acknowledgements.} We  thank J. Wahl for some
helpful conversations, and G. Faltings for a helpful communication. This work
was partially done while the
first author was visiting the  Max Planck Institut f\"ur Mathematik
(Bonn) in June, 1993 and the  Tata Institute of Fundamental Research (Bombay)
in
Dec., 1993; hospitality of which is gratefully acknowledged. The first author
lectured on the contents of this paper at the University of North Carolina
(Chapel Hill) in the fall of 1994,  at the Rutgers University (New Brunswick)
in Jan., 1995, and
at the University of Aarhus in May, 1995.  The first
author was supported by the NSF grant no. DMS-9203660.
\vskip4ex
   \head 1. Notation \endhead

 Let $G$ be a simple simply-connected connected complex affine
 algebraic group and let $C$ be a smooth irreducible
 projective curve  of genus $\geq 2$ over the field of complex numbers $\Bbb
 C$.  As in \cite{KNR, Theorem 3.4},
 let $\frak M$ be the moduli space of semistable
 principal $G$-bundles on $C$.  Also, fix a point $p\in
 C$ and recall the definition of the generalized flag
 variety $X=\Cal G/\Cal P$ (associated to the affine Kac-Moody
 group $\Cal G$ corresponding to the group $G$) from \cite{KNR,
 \S2.1}, its open subset $X^s$ and the morphism $\psi:X^s
 \to \frak M$ from \cite{loc. cit., Definition 6.1}.  Also,
 recall the
 notation $\widetilde{W},W,X_{\frak w}$ from  \cite{loc. cit.,\S2.1}.

For any ind-variety $Y$, by an {\it algebraic vector bundle  of rank } $r$ over
$Y$, we mean an ind-variety $E$ together with a morphism $\theta : E \to Y$
such that  (for any $n$) $E_n \to Y_n $ is an algebraic vector bundle over the
(finite dimensional) variety $Y_n$ of rank $r$ ,  where $\{ Y_n\}$  is the
filtration of $Y$ giving the ind-variety structure  and $E_n := \theta ^{-1}
(Y_n)$ . If $r=1$ , we call $E$ an {\it algebraic line bundle} over $Y$.

Let $E$ and $F$ be two algebraic vector bundles over  $Y$. Then a  morphism (of
ind-varieties) $\varphi : E \to F$ is called
a {\it bundle morphism} if the following diagram is commutative :
   $$
       \gather
       E
       \overset\varphi\to{\longrightarrow}
       F \\
       \searrow \qquad \swarrow\\
       X
       \endgather
       $$
and moreover $\varphi_{|E_n} : E_n \to F_n $ is a bundle morphism for all $n$.
In particular, we have the notion of isomorphism of vector bundles over  $Y$.

We define Pic $Y$ as the set of isomorphism classes of algebraic line bundles
on $Y$. It is clearly an abelian group under the tensor product of line
bundles.

For any set $Y$, $I_Y$ denotes the identity map of $Y$.
 \vskip4ex

   \head 2. Statement of the Main Theorems \endhead

 We follow the notation from \S1.
 \proclaim{(2.1) Lemma } The morphism $\psi : X^s \to
 \frak M$ induces an injective map
 $$
 \psi^\ast : \text{Pic}(\frak M) \longrightarrow
 \text{Pic}(X^s).
 $$
 \endproclaim

 \demo{Proof} Let $\frak L\in \text{Pic}(\frak M)
 $ be in the kernel of
 $\psi^\ast$, i.e., $\psi^\ast (\frak L)$ admits a
 nowhere-vanishing regular section $\sigma$ on the whole
 of $X^s$.  Fix $m\in \frak M$ and a trivialization for
 $\frak L_{\vert _m}$.  This canonically induces a
 trivialization for the bundle $\psi^\ast (\frak
 L)_{\vert \psi^{-1}(m)}$.  In particular, the section
 $\sigma_{\vert \psi^{-1}(m)}$ can be viewed as a (regular)
 map $\sigma_m:\psi^{-1}(m) \to \Bbb C^\ast $.  But
 $\psi^{-1}(m)$ is a certain union of $\Gamma$-orbits say
 $\psi^{-1}(m)= \underset{i\in I}\to{\cup}\Gamma x_i$, for
 $x_i\in X$ and moreover $\overline{\Gamma x_i}\cap
 \overline{\Gamma x_j}\neq \emptyset$, for any $i,j\in I$
 (cf. \cite{KNR, Proof of Proposition 6.4}).  Fixing
 $i\in I$, we get a regular map $\sigma_{m,i}: \Gamma \to
 \Bbb C^\ast $, defined as $\sigma_{m,i}(\gamma
 )=\sigma_m(\gamma x_i)$, for $\gamma\in \Gamma $.  Now
 by \cite{KNR,  Corollary 2.6}, $\sigma_{m,i}$ is
 a constant map for any $i\in I$, and hence $\sigma_m:
 \psi^{-1}(m)\to \Bbb C^\ast$ itself is a constant map.
 Thus the section $\sigma$ descends to a set theoretic
 section $\hat{\sigma}$ of the line bundle $\frak L$,
 which is regular by \cite{KNR, Proposition 4.1 and Lemma
 6.2}.  Of course, the section $\hat{\sigma}$ does not
 vanish anywhere on $\frak M$ (since $\sigma$ was chosen
 to be nowhere-vanishing on $X^s$).  This proves that
 $\frak L$ is a trivial line bundle on $\frak M$, thereby
 proving the lemma.\qed
 \enddemo

It is clear that for any ind-variety $Y$, we have a natural map $\alpha :$ Pic
$Y \to  \varprojlim_ {n}\,
 \text{Pic}\, (Y_n)\,.
$

 \proclaim{(2.2) Lemma}   Pic$~X  ~
\approx
\varprojlim_ {\frak w\in\widetilde{W}/W}\,
 \text{Pic}\, (X_{\frak w})\, \approx \Bbb Z ~.
 $
\endproclaim
 \demo{Proof}  We will freely follow the notation from \cite{KNR, \S2.3}. Since
 the
 line bundles $\frak L(d\chi_0)$ (for $d \in {\Bbb Z}$) are, by construction,
algebraic line bundles on $X$ and moreover for any ${\frak w} \geq {\frak s}_o$
,   ${\frak L}(\chi_0)_{\vert X_ {\frak w}} $ freely generates Pic ($X_{\frak
w}$), the surjectivity of the map $\alpha$ follows.  Now we come to the
injectivity of $\alpha$ :

Let $\frak L \in $ Ker $\alpha$.  Fix a non-zero vector $v_o$ in the fiber of
$\frak L$ over the base point $\frak e \in X $.  Then
${\frak L}_{\vert  X_{\frak w}} $ being a trivial line bundle on each $X_{\frak
w}$,
  we can choose a  nowhere-vanishing section $s_{\frak w}$
of  ${\frak L}_{\vert X_\frak w} $ such that  $s_{\frak w} (\frak e) =  v_o .$
We next show that for any $\frak v \geq \frak w ,  ~s_{\frak v_{\vert X_{\frak
w}}}=  s_{\frak w} :$  Clearly  $ s_{\frak v_{\vert X_{\frak w}}} =  f s_{\frak
w} ,$ for some algebraic function $f : X_{\frak w} \to
{\Bbb C}^\ast$ . But $X_{\frak w}$ being projective and irreducible,
$f$ is constant and in fact $f\equiv 1$ since $s_{\frak v} (\frak e) = s_{\frak
w} (\frak e) .$ This gives rise to a nowhere-vanishing section
$s$ of $\frak L$ on the whole of $X$ such that  $s_{\vert X_\frak w}=  s_{\frak
w} $ . From this it is easy to see that $\frak L$ is isomorphic with the
trivial line bundle on $X$. This proves that $\alpha$ is an isomorphism. Now
the second  isomorphism is proved in \cite{KNR, Proposition 2.3}.
\qed
\enddemo

 We state the following very crucial \lq lifting' result,
 the proof of which will be given in the next section.

 \proclaim{(2.3) Proposition} There exists a  map
 $\overline{\psi^\ast }: \text{Pic}(\frak M)\to
 \text{Pic}(X)$, making the following diagram
 commutative:
 $$
   \gather
   \text{{\rm Pic}}\, (\frak M)\\
               \vspace{1\jot}
   \sideset{\overline{\psi^\ast}}\and\to\swarrow
  \qquad\quad\quad
  \sideset\and{\psi^\ast }\to\searrow\\
            \vspace{2\jot}
  \text{{\rm Pic}}(X) \quad \underset{i^\ast}\to \rightarrow
 \quad
 \text{{\rm Pic}}(X^s) ,\endgather
  $$
 where $i^\ast$ is the canonical restriction  map.
  \endproclaim

 As an easy consequence of the above proposition, Lemmas (2.1) and (2.2), we
get
 the following  main result of this paper.

 \proclaim{(2.4) Theorem} For any smooth projective
 irreducible curve $C$ of genus $\geq 2$ and simple simply-connected connected
affine
 algebraic group $G$, the map $\overline{\psi^\ast }$ (as
 in the above proposition) is an injective group
 homomorphism.

 In particular, {\rm Pic}$\,(\frak M)\approx \Bbb Z$.
 \endproclaim

 \demo{Proof}  Injectivity of $\overline{\psi^\ast }$
 follows from the injectivity of $\psi^\ast $ (cf. Lemma
 2.1) and the commutativity of the diagram in Proposition
 (2.3).  By Proposition (2.3), Image $\psi^\ast \subset $ Image $i^\ast$ . But
since Pic $X  \approx  {\Bbb Z} $ (by Lemma 2.2), Image $i^\ast$ is either
finite or else Image $i^\ast \approx {\Bbb Z}$ . Now since $\frak M$ is a
projective variety of dim $> 0 $  (cf. \cite{R1, Theorem 4.9})
and $\psi^\ast$ is injective (Lemma 2.1), Image $i^\ast$ can not be finite, in
particular, $i^\ast$ is injective. Since $\psi^\ast$ and $i^\ast $ are group
 homomorphisms and $i^\ast $ is injective,
 we get that $\overline{\psi^\ast }$ is a group
 homomorphism. This proves the theorem. ~~~ \qed
 \enddemo

 \flushpar {\bf (2.5) Definition.} Let $n_{_{C,G}}>0$ be the
 least (positive) integer such that $\frak L(n_{_{C,G}}\chi_0)
 \in $ Image$\,\overline{\psi^\ast }$.  Then of
 course
 $$
 \text{Image}\,\overline{\psi^\ast}
 =\{\frak L(dn_{_{C,G}}\chi_0)\}_{d\in
 \Bbb Z}
 $$
\vskip1ex
 We will be concerned with determining the number $n_{_{C,G}}$ in \S4.
 \vskip1ex
 \flushpar {\bf (2.6) Remark.}  In the case when $G=
 SL(n,\Bbb C)$, it is a result of Drezet--Narasimhan
 \cite{DN} that Pic$\,(\frak M)\approx \Bbb Z$.
 \vskip1ex

 We recall the following  well known result.  (We include a proof since we did
not find it in the literature in this form.)

 \proclaim{(2.7) Lemma}  Let $Y$ be a Cohen--Macaulay
 projective variety and let $U\subset Y$ be an open
 subset such that codim$_Y(Y\backslash U)\geq 2$.  Now
 let $\Cal S_1$ and $\Cal S_2$ be two reflexive sheaves on $Y$
 such that $\Cal S_{1_{\vert_{U}}}\approx \Cal S_{2_{\vert_{U}}}$.  Then
the sheaf $\Cal S_1$ is isomorphic
 with $\Cal S_2$ on the whole of $Y$.
 \endproclaim

\demo{Proof\footnote{This proof is due to N. Mohan Kumar.}} We recall the
following two facts from Commutative Algebra.

\noindent
Fact 1: If $M,N$ are modules with depth $M, N >1$, and
$0\to M\to N\to K\to 0$  is an exact sequence, then depth $K>0$.

\noindent
Fact 2: If $M$ is reflexive, then for any localisation $M_{\frak p}$ of $M$ at
a prime
ideal $\frak p$, depth $M_{\frak p}>1$, unless the dimension of the local ring
itself is less than $2$ (i.e. $M$
satisfies the `Serre condition' $S_2$).

Let $i: U \hookrightarrow Y$ be the inclusion. Then from the above facts (and
the assumptions of the lemma), one can check that $i_*i^* \Cal S_j=\Cal S_j$
(for $j=1,2$). Thus any
homomorphism $i^*\Cal S_1\to i^* \Cal S_2$  on U gives rise to a homomorphism
$\Cal S_1\to \Cal S_2$ , i.e.,
Hom $(\Cal S_1,\Cal S_2)\to$  Hom $(i^*\Cal S_1, i^*\Cal S_2)$  is surjective.
Injectivity is clear
using reflexivity. This proves the lemma. \qed
\enddemo

 We come to the following second main result of this
 paper.

 \proclaim{(2.8) Theorem} The dualizing sheaf $\omega$ of the
 moduli space $\frak M$ is locally free.  Moreover,
 $\overline{\psi^\ast }(\omega)= \frak L(- 2g \chi_0
 )$, where $g$ is the dual Coxeter number of the Lie
 algebra $\frak{g}$ (cf. \cite{KNR, Remark 5.3}).

 In particular, $\frak M$ is a Gorenstein variety.
 Further, for any line bundle $\frak L$ on $\frak M$ such
 that $\overline{\psi^\ast }(\frak L)= \frak L(d\chi_0)
 $ for some $d>-2g$, $H^i(\frak M,\frak L)=0$, for
 all ~ $i>0$.
So, for any finite dimensional representation $V$ of $G$ ,
 $H^i(\frak M,\Theta(V))=0$, for
 all ~ $i>0$.
\endproclaim

 \demo{Proof}  Let $\frak M^o := \{E \in \frak M ; E \text{ is a stable }
G-\text{bundle and Aut} E = \text{centre of } G \} .$ Then  $\frak M^o $ is an
open subset  of the smooth locus of $\frak M$  and  for any $E\in \frak M^o$,
 the tangent space $T_E(\frak M^o)$ can be identified
 with $H^1(C,\text{{\rm ad}}\,E)$, where ad$\,E$ is the
 vector bundle on $C$ associated to the principal
 $G$-bundle $E$ via the adjoint representation of $G$ in
 its Lie algebra  $\frak{g}$.  Also, on the set of stable bundles in the moduli
space there are no identifications, i.e., if $E_1$ and $E_2$ are
two stable $G$-bundles on $C$ such that $E_1$ is $S$-equivalent to $E_2$, then
$E_1$ is isomorphic with $E_2$ (as follows from the definition of
$S$-equivalence, cf. \cite{KNR, \S3.3}).  Moreover, for any $E\in
 \frak M^o$, $H^0(C,\text{{\rm ad}}\,E)=0$.  In
 particular, the fiber of the canonical bundle of $\frak
 M^o$ at $E$ can be identified with
 $\wedge^{\text{top}}(H^1(C,\text{{\rm ad}}\,E)^\ast )$,
 where $\wedge^{\text{top}}$ is the top exterior power.
 This gives, from the definition of the determinant
 bundle and the $\Theta$-bundle (cf. \cite{KNR, \S3.8}),
 that $$\text{Det}\,(\text{ad})^\ast _{\vert_{\frak{M}^o}} =
 \Theta(\text{ad})^\ast _{\vert_{\frak{M}^o}}
 =\omega_{\vert_{\frak{M}^o}} . $$ But
 $\Theta(\text{ad})^\ast $ is a line bundle on the whole
 of $\frak M$, in particular, it is a reflexive sheaf on
 $\frak M$ by \cite{H, Exercise 5.1, p. 123}.  Since the
 dualizing sheaf $\omega$ of a normal variety is always
 reflexive;  the moduli space $\frak M$ is
 Cohen--Macaulay and normal  (cf. \cite{R1, Theorem 4.9});  and  codim$_{\frak
 M} (\frak M \backslash \frak M^o)\geq 2$ (unless  the curve $C$ is of genus
$2$ and $G=SL(2)$) ( cf. \cite{F, Theorem II.6} ); we obtain
 from Lemma (2.7):
 $$
 \omega \approx \Theta(\text{{\rm ad}})^\ast ,\text{ on the
 whole of }\frak M.
 \tag1$$
 (In the case of $G=SL(2)$ the validity of (1)  is well known.) This of course
gives that $\frak M$ is a Gorenstein
 variety (by definition).  Now the assertion that
 $\overline{\psi^\ast }(\omega)= \frak L(-2g\chi_0)$
 follows from \cite{KNR, Theorem 5.4 and Lemma 5.2}.

 Finally we come to the proof of cohomology vanishing:
 By Serre duality \cite{H, Corollary 7.7, Chap. III},
 $$\align
 H^i (\frak M,\frak L)^\ast  &\approx H^{n-i} (\frak
 M,\frak L^\ast \otimes \omega) \\
 &= H^{n-i} (\frak M, \frak L^\ast \otimes \Theta
 (\text{ad})^\ast ), \text{ by (1).} \tag2\\
 \endalign$$
 But $\overline{\psi^\ast } (\frak L^\ast \otimes \Theta
 (\text{ad})^\ast ) = \frak L((-d-2g) \chi_0)$.  Now
 since Pic$\,(\frak M)\approx \Bbb Z$  (by Theorem 2.4),
 we get that the line bundle $\frak L \otimes \Theta
 (\text{ad} )$ is ample on $\frak M$ (by assumption $d > -2g$).

The moduli space $\frak M$ has rational singularities, as follows from \cite{R,
Proof of Theorem 4.9} and a result of Boutot \cite{Bo}. Now the
 vanishing of $H^i(\frak M,\frak L)$ (for $i>0$) follows
 from (2) and a result of Grauert-Riemenschneider \cite{GR}. So the proof of
the theorem is complete in view of \cite{KNR, Theorem 5.4}.\qed
\enddemo

 \proclaim{(2.9) Corollary}  For any finite dimensional
 representation $V$ of $G$,
 $$
 \Cal X(\frak M,\Theta(V)) = \dim H^0(\frak M,\Theta(V)),
 $$
 where $\Cal X$ is the Euler-Poincare characteristic:
 $$
 \Cal X(\frak M,\Theta(V)) =
 \sum_i (-1)^i \dim H^i(\frak M,\Theta(V)).
 $$
 \endproclaim

\vskip3ex
 \head 3.  Extension of Line Bundles - Proof of
 Proposition (2.3) \endhead

 \subhead (3.1) \endsubhead  Recall the definition of
 the map $\varphi :\Cal G \to \Cal X_0$  from \cite{KNR,
 \S1} (where $\Cal
 X_0$ denotes the set of isomorphism classes of
 principal $G$-bundles on $C$ which are algebraically
 trivial restricted to $C^\ast $).  Fix an embedding $G\hookrightarrow SL(n) $,
for some
 $n$.  In particular, any principal $G$-bundle $E$ on
 $C$ gives rise to a vector bundle $\overline{E}$ of
 rank $n$ on $C$ (associated to the standard
 representation of $SL(n)$). For any integer $d\geq 1$,
 define
 $$
 X_d= \{ g\Cal P\in X:  H^1(C,\overline{\varphi (g)}
 \otimes \Cal O(-x+dp)) =0 \text{~for~all~} x\in C\},
 $$
 where $p\in C$ is the fixed base point.  Then
 $$
 X_1\subset X_2\subset \cdots .$$

\proclaim{(3.2) Lemma}  Each $X_d$ is open in $X$. Moreover
$X^s \subset X_{2g}$,  where
  $X^s=\{
 g\Cal P\in X: \varphi (g) \text{ is a semistable
 $G$-bundle}\}$,   and $g$ is the genus
 of the curve $C$.

\endproclaim
\demo{Proof}   It suffices to  prove that $X_d \cap X_{\frak w}$ is open in
$X_{\frak w}$ , for each $\frak w \in \widetilde{W}/W$ :

Recall the definition of the family of $G$-bundles
 $\Cal U \to C\times X$ from \cite{KNR, Proposition
 2.8}. Consider the restriction $\Cal U_{\frak w} $ of the $G$-bundle
 $\Cal U \to C\times X$ to $C \times  X_{\frak w} $ and let $\overline{\Cal
U_{\frak w}}$ be the associated rank-$n$ vector bundle (corresponding to the
embedding $G \hookrightarrow \text{SL} (n)$).
Define a vector bundle  $\widetilde{\Cal U_{\frak w}}$ on $C\times C \times
X_{\frak w} $ such that  $\widetilde{\Cal U_{\frak w}}_{\vert x\times C \times
X_{\frak w} }= \Cal O (-x+dp)  \otimes \overline{\Cal U_{\frak w}}$ ~ for each
$x \in C$ ;  and  let $\pi : C\times C \times X_{\frak w} \to  C \times
X_{\frak w} $ be the projection on the two extreme  factors. Applying the upper
semi-continuity theorem  \cite{H,Chapter III, \S12} to the morphism $\pi$ and
the locally free sheaf  $\widetilde{\Cal U_{\frak w}}$ on $C\times C \times
X_{\frak w} $, we get that the set
$$ S := \{(x, g\Cal P) :  H^1(C,\overline{\varphi (g)}
 \otimes \Cal O(-x+dp)) \neq 0 \} $$
is a closed subset of  $C \times X_{\frak w} $. In particular, $\pi_2 (S)$
is a closed subset of  $ X_{\frak w} $, where $\pi_2 : C \times X_{\frak w} \to
 X_{\frak w} $ is the projection on the second factor. It is easy to see that
$X_d
\cap X_{\frak w} =  X_{\frak w} \setminus \pi_2 (S) $ . This proves that $X_d$
is open in $X$.

For $g\Cal P \in X^s ,  ~ \overline{\varphi (g)} $ is a semistable vector
bundle, and hence the dual vector bundle $\overline{\varphi (g)}^\ast $
is also semistable. Now, by the Serre duality,
$$ H^1(C,\overline{\varphi (g)}
 \otimes \Cal O(-x+dp)) \approx  H^0(C,\overline{\varphi (g)}^\ast
 \otimes \Cal O(x-dp) \otimes K) ^\ast ~.$$
Since  $\overline{\varphi (g)}^\ast  $ is semistable, $H^0(C,\overline{\varphi
(g)}^\ast
 \otimes \Cal O(x-dp) \otimes K) \neq 0 $ implies that  $d-1-$ deg $K \leq 0$ .
In particular, if $d \geq 2 + \text{deg} ~K  ,$ then $g \Cal P \in X_d$ . This
proves  the lemma since deg $K = 2g-2$. \qed
\enddemo

We have
$$\underset {d\geq 1}\to{\cup} X_d = X ~;
 $$
in fact each Schubert variety $X_{\frak w}$ is contained in some large enough
$X_d$ ($d$ of course
depending upon $\frak w$). This follows by the upper semi-continuity theorem
(using an argument similar to the one used in the proof of the above lemma).

\vskip2ex
\noindent
 {\bf  (3.3)}   Fix any $d\geq 2g$. For all $m\geq d$ and $g\Cal P\in X_d$, we
have
 \vskip2ex
\roster
 \item $H^1(C,\overline{\varphi (g)} \otimes \Cal
 O(mp))=0$,  and
 \item $H^0(C,\overline{\varphi (g)} \otimes \Cal
 O(mp))$ generates the vector bundle
 $\overline{\varphi (g)} \otimes \Cal O(mp)$  at every
 point of $ C$.
 \endroster

 Let $q_d:=\dim H^0(C, \overline{\varphi (g)} \otimes
 \Cal O(dp))$.  Then by Riemann-Roch theorem, $q_d=
 n(d+1-g)$.  Denote
 by $\pi _d: \Cal F_d\to X_d$ the GL($q_d$)-bundle such
 that for $g\Cal P\in X_d$, $\pi _d^{-1}(g\Cal P)$ is the set of
 all the frames of the vector space
 $H^0(C,\overline{\varphi (g)} \otimes \Cal O(dp))$.
 We call $\Cal F_d$ the frame bundle associated to the
 family $\Cal U _{\vert_{X_d}}$ parametrized by $X_d$.
 Similarly,  define the frame
 bundle $\pi _{d+1}: \Cal F_{d+1} \to X_{d+1}$.
 Consider the parabolic subgroup $P=\{\theta\in$
 GL$(q_{d+1}): \theta\Bbb C^{q_d}= \Bbb C^{q_d}\}$ of
 GL$(q_{d+1})$, where (for definiteness) $\Bbb
 C^{q_d}\hookrightarrow \Bbb C^{q_{d+1}}$ is sitting in the
 first $q_d$ coordinates.  We define the principal
 $P$-subbundle $Q_d$ of ${\Cal F_{d+1}}_{\vert_{X_d}}$ by
 $$\multline
 Q_d =
 \underset {g\Cal P\in X_d}\to{\cup}
 \{ s=(s_1,\dots ,s_{q_{d+1}}) \text{ a frame of }
 H^0 (C, \overline{\varphi (g)} \otimes \Cal
 O((d+1)p)) \text{ such that }\\
 (s_1,\dots ,s_{q_{d}}) \text{ is a frame of }
 H^0 (C, \overline{\varphi (g)} \otimes \Cal
 O(dp)) \}.
 \endmultline$$
(Observe that   $H^0 (C, \overline{\varphi (g)} \otimes \Cal
 O(dp))$ sits canonically inside  $H^0 (C, \overline{\varphi (g)} \otimes \Cal
 O((d+1)p))$ induced from the embedding  $ \overline{\varphi (g)} \otimes \Cal
 O(dp) \hookrightarrow \overline{\varphi (g)} \otimes \Cal
 O((d+1)p)$.)
 Then we have the following commutative diagram :
 $$ \gather
    \Cal F_d
 \quad \quad
   \overset{\beta_d}\to\twoheadleftarrow
 \quad
  Q_d \quad \hookrightarrow
 \quad
    \Cal F_{d+1}\\
  \vspace{2\jot}
    \sideset{\pi _d}\and\to\downarrow
 \qquad\qquad\qquad\qquad\qquad
  \sideset\and{\pi _{d+1}}\to\downarrow\\
 \vspace{1\jot}
    X_d \qquad\qquad \hookrightarrow \qquad X_{d+1},
   \endgather
  $$
 where $\beta _d$ takes any $s=(s_1,\dots
 ,s_{q_{d+1}})\in Q_d$ to the frame $(s_1,\dots
 ,s_{q_d})$ of $H^0(C, \overline{\varphi (g)} \otimes
 \Cal O(d p))$.  It is clear that $\beta _d$ is a
 principal $U$-bundle, where $U:= \{ \theta\in
 \text{{\rm GL}}(q_{d+1}) : \theta_{\vert_{\Bbb C^{q_d}}}=
 I\,\}\subset P$.  Then clearly $U$ is a
 normal subgroup of $P$.

  As in \cite{KNR, \S7.8}, we have
 an irreducible smooth quasi-projective variety $R_d$
 with an action of GL$(q_d)$, a family $\Cal W_d$ of
 $G$-bundles on $C$ parametrized by $R_d$ and a lift
 of the GL$(q_d)$-action to $\Cal W_d$ (as bundle
 automorphisms) such that there exists a
 GL$(q_d)$-equivariant morphism  $\varphi _d:\Cal
 F_d\to R_d$ with the property that the families
 $\pi^\ast _d(\Cal U_{\vert_{X_d}})$ and $\varphi
 ^\ast _d(\Cal W_d)$ are isomorphic.  Moreover, let
 $R^s_d= \{ x\in R_d: \Cal W_d(x):= {\Cal
 W_d}_{\vert_{C\times x}} \text{ is a semistable
 $G$-bundle}\}$ be the
 GL$(q_d)$-invariant open subset
 of $R_d$.  Then the canonical map $\theta_d:R^s_d\to
 \frak M$ is surjective.  Moreover, $\theta_d$ is
 GL$(q_d)$-equivariant  with respect to the trivial
 action of GL$(q_d)$ on the moduli space $\frak M$ (of
 semistable $G$-bundles on $C$). We recall the construction of $R_d$ for its
use in the sequel \cite{R1, \S\S 3.8,3.13.3}:

Let $R_d^o$ be the set of locally free quotients $E$ of ${\Bbb
C}^{q_d}\otimes_{\Bbb C} \Cal O_C$ of rank $n$ and degree $nd$ such that the
canonical map ${\Bbb C}^{q_d} \approx H^0({\Bbb C}^{q_d}\otimes_{\Bbb C} \Cal
O_C) \to H^0(E) $
is an isomorphism. Then $R_d^o$ supports a tautological family
$\widehat{\Cal W_d}^o $ of rank-$n$ vector bundles on $C$. Set
$\Cal W_d^o = \widehat{\Cal W_d}^o \otimes_{\Cal O_{C \times R_d^o}}
\Cal O_C(-d p) $ . Now let
$$R_d = \{(x,\sigma): x\in R_d^o \text{~and~} \sigma \text{~is ~ a ~ reduction
{}~ of ~ the ~ structure ~ group ~ of } {\Cal W_d^o}_{\vert C \times x}
\text{~to~} G \}.$$
Then clearly $R_d$ supports a family $\Cal W_d$ of $G$-bundles on $C$ and
moreover GL($q_d$) acts canonically on   $\Cal W_d$ via
its action on ${\Bbb C}^{q_d}$.

Using $H^1(C,E)=0$, one proves that $R_d$ is smooth and that the infinitesimal
deformation map $T_t (R_d) \to H^1(C, \text{Ad}~{\Cal W_d}_{\vert C\times t})$
is surjective, where $T_t (R_d)$ is the tangent space at $t$ to $R_d$.

 \proclaim{(3.4) Proposition} For any $d\geq 2g$,    the codimension
 of $R_d\backslash R_d^s$ in $R_d$ is at least 2, where $R_d$  is explicitly
constructed as above.
\endproclaim

 To prove the above proposition, we need the notion of the
 canonical reduction (or  filtration) of a principal
 $G$-bundle on $C$.
 We choose a Borel subgroup $B$ of $G$ and a maximal torus
 $T\subset B$.  By a {\it standard parabolic subgroup}  we mean  a
 parabolic subgroup $P$ containing $B$.  The following result is due to
Ramanathan \cite{R2, Proposition 1} (see also \cite{Be}).

 \proclaim{ (3.5) Theorem} Let $E$ be a principal $G$-bundle on
 $C$.  Then there exists a unique standard parabolic subgroup $P$ of
 $G$ and a unique reduction $E_P$ of $E$ to the
  subgroup  $P$ such that the following
 conditions hold:
 \roster
 \item If $U$ is the unipotent radical of $P$, then the
 $P/U$-bundle $E_{P/U}$ obtained from $E_P$ by extension
 of the structure group via $P\to P/U$ is semi-stable.
 (Observe that $P/U$ is reductive.)
 \item For any non-trivial character $\chi $ of $P$ which is a
 non-negative linear combination  of simple roots of $B$, the line bundle
 on $C$ associated to $E_P$ by $\chi $ has
 strictly positive degree.
 \endroster \endproclaim

The unique reduction $E_P$ of $E$ as above is called the {\it canonical
reduction}.
  \proclaim{(3.6) Lemma}  Let $E_P$ be the canonical reduction
 of a principal $G$-bundle $E$ on $C$.  Let $\frak g$ and $\frak
 p$ be the Lie algebras of $G$ and $P$ respectively.
 Denote by $E_{\frak s}$ the vector bundle associated to $E_P$ by
 the natural representation of $P$ on the vector space
 $\frak s :=\frak g/\frak p$.  Then we have
 $$
 H^0(C,E_{\frak s})=0.
 $$
 \endproclaim

 {\it Proof}.  We may assume that $P\neq G$.  Let
 $0=V_0\subset V_1 \dots  \subset V_k= \frak s$ be a
 filtration of $\frak s$ by $P$-submodules $V_i$ such that for any $1\leq i
\leq k$ , the $P$-module $W_i:=V_i/V_{i-1}$ is irreducible. In particular,  $U$
acts trivially on $W_i$ (cf. \cite{Ku, Lemma 1}).  If $\Cal V _i$ is the vector
bundle on $C$
 associated to $E_P$ by the representation of $P$ on
 $V_i$, then $E_{\frak s}$ is filtered by the subbundles $\Cal V
 _i$.  We now show that $H^0(C,\Cal W _i)=0$ for all $1\leq i \leq k$, where
$\Cal W _i := \Cal V _i / \Cal V _{i-1}$.   This will of course prove
 the lemma.

 Since the action of $U$ on $W_i$ is trivial,  we
 obtain an (irreducible)  representation of the reductive group $P/U$ on
 $W_i$.
 Since $E_{P/U}$ is semi-stable, the vector bundles $\Cal
 W_{i}$ are semi-stable, and hence  it is
 sufficient to show that $\deg (\Cal W_{i})<0$.  Now
 the weights of $T$ on $\frak s$ are
of the form
 $\sum c_\alpha\alpha$ with $c_\alpha\leq 0$ and
 $c_\alpha<0$ for at least one $\alpha\notin I$, where
 $I$ is the subset of simple roots $ \Pi=\{\alpha\}$ defining
 the parabolic subgroup $P$ (i.e. $I$ is the set of simple roots for $P/U$). It
follows from this that the
 character of $P$ defined by the determinant of the
 representation of $P$ on $ W_{i}$ is non-trivial
 and is a non-positive linear combination of $\{\alpha \}_{\alpha \notin I}$.
By
 Condition 2) of Theorem (3.5),  we see that $\deg(\Cal
 W_{i})<0$.  This completes the proof of the lemma.
 \qed
 \vskip1ex

 Let  $P$ be  a standard parabolic subgroup of $G$ and
 $E_P$ be a reduction of the $G$-bundle $E$ to $P$.  For any  character $\chi$
of $P$, denote by $E_{P,\chi}$ the
 line bundle on $C$ associated to $E_P$ by $\chi$.  Let
 $X(P)$ (resp. $X(T)$) denote the character group of $P$
 (resp. $T$).  Then $X(T)= \oplus_{\alpha \in \Pi} \Bbb Z \omega_\alpha$ ,
where $\omega_\alpha$ is the fundamental weight defined by $\omega_\alpha
(\beta^\vee) = \delta_{\alpha,\beta}$ , for any simple coroot $\beta^\vee$.
Moreover (since $G$ is simply-connected)
 $X(P)= \oplus_{\alpha \notin I} \Bbb Z \omega_\alpha$ . The map $\chi \mapsto
\deg (E_{P,\chi})$
 defines an element of $\text{Hom}_{\Bbb Z}(X(P),\Bbb Z)$,  which in turn
 can be lifted  to the element $\mu $  of $\text{Hom}_{\Bbb Z} (X(T),\Bbb
 Z)$ defined by $\mu(\omega_\alpha) =$ deg $(E_{P,\omega_\alpha})$ if $\alpha
\notin I$ and $\mu(\omega_\alpha) = 0$ if $\alpha \in I$.  We call $\mu $ the
{\it type}  of the reduction $E_P$.

 Using the above  lemma, one can prove the following
 proposition; the proof being similar to that of \cite{PV, Theorem  4, p. 90}.

 \proclaim{(3.7) Proposition}  Let $\Cal W$ be a family of
 $G$-bundles on $C$ parametrized by a smooth variety $S$.
 Assume that at each point $t\in S$ the infinitesimal
 deformation map
 $$
 T_t(S) \rightarrow H^1 (C, \text{Ad}(\Cal W_t))
 $$
 is surjective, where $\Cal W_t=\Cal W_{\vert_{C\times t}}$
 and $T_t(S)$ is the tangent space at $t$ to $S$.  For
 $\mu\in \text{Hom}(X(T), \Bbb Q)$, let $S_{\mu}$ be the subset of
 $S$ consisting of those points $t \in S$ such that the canonical
 reduction of $\Cal W_t$ is of type $\mu $.  Then $S_\mu$ is
 non-empty only for finitely many $\mu $.  Moreover,
 $S_\mu$ is locally closed  and smooth, and the normal space at $t\in
 S_\mu$ is given by $H^1(C, \Cal W_{t,\frak s})$, where $\Cal W_{t,\frak s}$
 is the vector bundle associated to the canonical
 reduction $\Cal W_{t,P}$ by the representation of $P$ on
 $\frak s :=\frak g/\frak p$.
 \endproclaim

 \noindent
{\bf (3.8)} {\it Proof of Proposition (3.4).}
 The family $\Cal W=\Cal W_d$ parametrized by $R_d$
 satisfies the hypothesis of the above proposition (3.7).  So it
 suffices to prove that for  $t\in R_d\backslash R_d^{s}$,
 we have $\dim H^1(C,\Cal W_{t,\frak s}) \geq 2$. By  Lemma (3.6),  $H^0(C,\Cal
W_{t,\frak s})=0$  and we have by Riemann-Roch theorem,
 $$
 \dim H^1(C,\Cal W_{t,\frak s})= - \deg \Cal W_{t,\frak s} +
 \dim(\frak s) (g-1) \, , \tag 1
 $$
 where recall that $g$ is the genus of $C$. Further,  since $t\in R_d\backslash
R_d^s$ ,  we have $\frak
 g\neq \frak p$.  The weight of  $\wedge^{\text{top}}\frak s = -2
 \rho_o $ ,  where $\rho_o := \sum_{\alpha \notin I} \omega_\alpha $ .   Write
 $$
 2\rho_o = \sum_{\alpha\in  \Pi} n_\alpha \alpha ~,
 \tag 2 $$
 for some non-negative integers  $n_\alpha$ .
 So
 $$
 \deg \Cal W_{t,\frak s} = -2 \deg(E_{P,
 \rho_o})~. \tag 3
 $$
 (Observe that $\rho_o \in X(P)$  and moreover it is a non-trivial character.)
By Condition (2) of Theorem (3.5) and (2), $\deg(E_{P,\rho_o})\geq 1$, and
hence by (3)
 $$
 \deg \Cal W_{t,\frak s} \leq  -2.
 $$
 This gives  (using 1) that $\dim H^1(C,\Cal W_{t,\frak s})\geq 2$, proving
Proposition (3.4). \qed

\proclaim{(3.9) Lemma} Let $H$ be an affine algebraic group
acting algebraically on a smooth variety $Y$ and let $U$ be a
$H$-stable open subset such that ${\text codim}_Y(Y\backslash U)
\geq 2.$
 Then the canonical restriction map {\rm
 Pic}$^H(Y)\to\,${\rm Pic}$^H( U)$ is an
 isomorphism, where {\rm Pic}$^H(Y)$ denotes the set
 of isomorphism classes of $H$-equivariant line bundles
 on $Y$.\endproclaim

 \demo{Proof}  Let $\Cal L$ be an $H$-equivariant line
 bundle on $ U $.  Since $Y$ is smooth and
 codim$_Y(Y\backslash  U ) \geq 2$, $\Cal L$ extends
 uniquely to a line bundle $\tilde{\Cal L}$ on $Y$.
 We show that $\tilde{\Cal L}$ is $H$-equivariant:

Fix $h\in H$ and an open subset $V\subset Y$ such that  $\tilde{\Cal
L}_{\vert_{V}}$ is a trivial line
 bundle.  In particular, the line bundle $\tilde{\Cal
 L}_{\vert_{hV}}$  also is trivial (since by the
 $H$-equivariance of $\Cal L$, $\tilde{\Cal
 L}_{\vert_{h( U \cap V)}}$ is trivial and
 moreover codim$_V(V\backslash  U )\geq
 2$).  Take a nowhere-vanishing section $s_1$ of
 $\tilde{\Cal L}_{\vert_{V}}$ and $s_2$ of
 $\tilde{\Cal L}_{\vert_{hV}}$.  Now for any $x\in
  U \cap V$, $f_h(x)s_2(hx)= h(s_1(x))$, for some
 (unique) $f_h(x)\in \Bbb C^\ast $.  Clearly the map
 $ U \cap V\to \Bbb C^\ast $, taking $x\mapsto
 f_h(x)$ is a regular map, which extends to a regular
 map $\tilde{f_h}: V\to \Bbb C^\ast $ (since
 codim$_V(V\backslash  U )\geq 2$).  Define an
 action of $h$ on  $\tilde{\Cal L}_{\vert_{V}}$ by
 $$
 h(s_1(x)) = \tilde{f_h} (x)s_2(hx),\qquad \text{for
 all }x\in V.
 $$

 By the uniqueness of extension, this action of $h$
on $\tilde{\Cal L}_{\vert_{V}}$ patches-up to give an action of $h$ on the
whole of
 $\tilde{\Cal L}$. Further, as can be easily seen, this is a regular action
of $H$ on   $\tilde{\Cal L}$.

 The injectivity of Pic$^H(Y)\to \,$Pic$^H( U )$
 is easy to see:  An $H$-equivariant section, which
 does not vanish anywhere on $ U $, extends to a
 nowhere-vanishing section on $Y$ (and by uniqueness
 of extension it is $H$-equivariant).  \qed
 \enddemo

 \vskip2ex

 \subhead (3.10) Lifting of line bundles from  $\frak M$
 to $X_d$ \endsubhead
 Take any $d \geq 2g$. Let $\frak L$ be a line bundle on $\frak M$.  Pull back
the
 line bundle $\frak L$ via the  GL$(q_d)$-equivariant
 morphism $\theta_d:R^s_d\to\frak M$ to get a
 GL$(q_d)$-equivariant line bundle $\theta_d^\ast(\frak
 L)$ on $R^s_d$ (cf. \S3.3).  By the above Lemma (3.9) and
 Proposition (3.4), $\theta^\ast_d(\frak L)$ extends to a
 GL$(q_d)$-equivariant line bundle $\widehat{\theta^\ast
 _d(\frak L)}$ on $R_d$.  Consider the diagram, where
 all the maps are GL$(q_d)$-equivariant (the map $i_d$ is the inclusion,
$\varphi_d$ and $\pi_d$ are as in \S3.3, and GL$(q_d)$ acts trivially on
$X_d$):

$$ \gathered
\CD
  \Cal F_d  @>\varphi_d>>  R_d\\ @V\pi_dVV  @.\\  X_d
\endCD
\CD @<i_d\ \ <<  R^s_d\\  @.  @VV\theta_dV\\  @. \frak M
\endCD
\endgathered
$$

 Now $\varphi ^\ast _d(\widehat{\theta ^\ast _d(\frak
 L)})$ being a GL$(q_d)$-equivariant line bundle and
 $\pi_d$ is a principal GL$(q_d)$-bundle, this descends
 to give a line bundle (denoted) $\frak L_d$ on $X_d$
 (cf. \cite{Kr, Proposition 6.4}).

 \proclaim{(3.11) Lemma} For any line bundle $\frak L$ on
 $\frak M$ and $d \geq 2g$
 $$
 {\frak L_{d+1}}_{\vert_{X_d}} \approx \frak L_d ~,
 $$
and ${\frak L_d}_{\vert_{X^s}} \approx \psi^\ast(\frak L)$, where $\psi:X^s
 \to \frak M$ is the morphism as in \S1 (cf. Lemma 3.2).
 \endproclaim

 \demo{Proof} We will freely use the notation from \S3.3.
 Let $X_{\frak w}$ be a fixed Schubert variety, and
 denote   the (reduced) variety  $X_{\frak w}\cap X_d$ by $Y=Y_{d,\frak w}$.
 Then $Y^s :=Y\cap X^s$ is an open non-empty   (irreducible) subvariety of
$X_{\frak w}$.  We denote by
 $\Cal F_{d,Y}$, $\Cal F_{d+1,Y}$ and $Q_{d,Y}$ the
 restrictions of $\Cal F_d,~\Cal F_{d+1}$ and $Q_d$ to $Y$, where $Q_d$ is the
$P$-subbundle of ${\Cal F_{d+1}}_{\vert X_d}$ as in \S3.3.  We
  show that ${\frak L_d}_{\vert Y} \approx {\frak
 L_{d+1}}_{\vert Y}$ and  ${\frak L_d}_{\vert Y^s}\approx
 \psi^\ast(\frak L)_{\vert Y^s}$.  This will of course prove the lemma.

 We first show that
 $$
 {\frak L_d}_{\vert Y^s} \approx \psi^\ast(\frak L)_{\vert Y^s}~: \tag 1
 $$
 From the commutativity of the diagram (where $\Cal F^s _{d,Y}:=
\pi_d^{-1}(Y^s)$, and $\pi_d, \varphi_d,$ and $\psi$ are the corresponding
restriction maps, which we denote by the same symbols)
   $$
       \matrix
       && \Cal F^s _{d,Y}\\ \vspace{1\jot}
       &\sideset\pi_d\and\to\swarrow  &&\sideset\and\varphi_{d}\to\searrow \\
   \vspace{1\jot}
       &Y^s &&R^s_{d}\\ \vspace{1\jot}
       &\sideset\psi \and\to\searrow
  &&\sideset\and\theta_{d}\to\swarrow \\ \vspace{1\jot}
       &&\frak M \\
       \endmatrix \tag $D_1$
       $$
 we see that the GL$(q_d)$-linearizations on $\pi
 ^\ast_d(\psi^\ast\frak L)$ and
 $\varphi^\ast_d(\theta^\ast_d\frak L)$ are the same. This
 shows that ${\frak L_d}_{\vert Y^s}\approx \psi^\ast(\frak L)_{\vert Y^s}$
(since $\pi_d$ is a principal GL $(q_d)$-bundle).

 If $H$ is an affine algebraic group and $\Cal H$ an
 $H$-linearized line bundle on a principal $H$-bundle, we
 denote by $\Cal H^H$ the line bundle on the base space
 (of the $H$-bundle) obtained by descending $\Cal H$.

 Let $\widetilde{\Cal W}^o_d$ be the vector bundle on $C\times R_d$
 which is the pull-back of $\widehat{\Cal W}^o_d$ by the map
 $I_C\times \beta:C\times R_d\to C\times R^o_d$, where
 $\beta:R_d\to R^o_d$ is the canonical map.

 Let $\pi''_d :\Cal F''_{d} \to R_d$ (resp. $\pi'_d : \Cal F'_d \to R_d$) be
the frame bundle of
 the vector bundle
 $(p_{R_d})_\ast (\widetilde{\Cal W}^o_d\otimes \Cal O (p))$ (resp.
 $(p_{R_d})_\ast (\widetilde{\Cal W}^o_d)$), where $p_{R_d}: C\times R_d \to
R_d$ is the projection on the second factor.  Just as in \S3.3, the inclusion
 $$
 (p_{R_d})_\ast  (\widetilde{\Cal W}_d^o) \hookrightarrow
 (p_{R_d})_\ast  (\widetilde{\Cal W}_d^o \otimes \Cal O (p))
 $$
 defines a $P$-subbundle $Q'_d \subset \Cal F''_d$ on $R_d$ and a morphism
$\beta'_d : Q'_d \to \Cal F'_d$ . Further, analogous to the map $\varphi_d :
\Cal F_d \to R_d $ (cf. \S3.10), there is a GL$(q_{d+1})$-equivariant morphism
$\varphi'_d : \Cal F''_d \to R_{d+1}$ . Thus we have the
 diagram:
  $$
   \matrix\format\c&&\c\\
 &&Q'_d\\
\vspace{1\jot}
 \sideset\beta'_d\and\to\swarrow   &&&\searrow\\ \vspace{1\jot}
 \Cal F'_d &&& \Cal F''_{d}\\ \vspace{1\jot}
 \pi'_d @ VVV &@VV\varphi'_{d}V\\ \vspace{1\jot}
 R_d &&&R_{d+1}~.
   \endmatrix \tag $D_2$
  $$
 (Observe that $\beta'_d$ is a principal $U$-bundle, $\pi'_d$ is a principal
GL$(q_d)$-bundle and $\pi''_d$ is a principal GL$(q_{d+1})$-bundle.)
Considering the commutative diagram (where $\Cal F^{\prime\prime \,s}_d :=
\pi_d^{''-1} (R_d^s)$)
       $$
       \matrix
       && \Cal F^{\prime\prime \,s}_{d}\\ \vspace{1\jot}
       &\swarrow  &&\sideset\and\varphi'_{d}\to\searrow \\
   \vspace{1\jot}
       &R^s_d &&R^s_{d+1}\\ \vspace{1\jot}
       &\sideset\theta_d\and\to\searrow
  &&\sideset\and\theta_{d+1}\to\swarrow \\ \vspace{1\jot}
       &&\frak M\\
       \endmatrix \tag $D_3$
       $$
 we see, as above, that
 $$
 ( \varphi^{\prime\,\ast} _{d} \theta^\ast_{d+1} \frak
 L)^{\text{GL}(q_{d+1})} \approx \theta^\ast_d(\frak L) .
 $$
 Since codim$_{R_d}(R_d\backslash R^s_d)\geq 2$ and $R_d$ is
 smooth, we have
 $$
 \widehat{\theta^\ast_d \frak L } \approx
 (\varphi_{d}^{\prime \ast}
 (\widehat{\theta^\ast_{d+1} \frak L}))^{\text{GL}(q_{d+1})} .
 $$
 Now
 $$
 \align
  (\varphi^{\prime\,\ast}_{d} (
 &\widehat{\theta^\ast_{d+1} \frak L}
 ))^{\text{GL}(q_{d+1})} \\
 &\approx (\gamma^\ast_{d}
 (\widehat{\theta^\ast_{d+1} \frak L}))^P \\
 &\approx ((\gamma^\ast_{d}(\widehat{\theta^\ast_{d+1} \frak
 L}))^U)^{\text{GL}(q_d)}  \\
 &\approx \sigma^\ast ( (\gamma^\ast_{d}
 (\widehat{\theta^\ast_{d+1} \frak L}))^U ) ~,
 \endalign
 $$
 where $ \gamma_d : Q'_d \to R_{d+1}$ is the restriction of $\varphi'_d$ to
$Q'_d$ and $\sigma : R_d\to \Cal F'_d$ is the canonical section, given by the
 isomorphism
 $$
 \Bbb C^{q_d} = H^0 (C,\Bbb C^{q_d}\otimes \Cal O_C)
 \rightarrow H^0(C,\widetilde{\Cal W}^o_d\vert_{C\times t})
 $$
 for $t\in R_d$ .
 Thus
 $$
 \widehat{\theta^\ast_{d} \frak L} \approx
 \sigma^\ast ((\gamma^\ast_{d} (\widehat{\theta_{d+1}^\ast \frak
 L}))^U ) .
 \tag 2 $$

Consider the following  commutative diagram
 $$\matrix\format\c&&\quad\c \\
 Q_{d,Y} &\overset\alpha\to\longrightarrow &Q'_d  &\hookrightarrow &\Cal
F''_{d} \\
 \vspace{1\jot} \beta_d @VVV @VVV @VV\varphi'_dV\\ \vspace{1\jot}
 \Cal F_{d,Y} &\overset\delta\to\longrightarrow &\Cal F'_d &&R_{d+1}\\
\vspace{1\jot}
 \pi_d @VVV  @VVV\\ \vspace{1\jot}
 Y  &&R_d
 \endmatrix \tag $D_4$ $$
 where  $\delta :=\sigma \circ \varphi _d$, and
the map $\alpha$ is  defined as follows:
 Let $g\Cal P\in Y$ and $s=(s_1, \dots , s_{q_d}, \dots ,
 s_{q_{d+1}})$ be a frame of $H^0(C,\overline{\varphi (g)}\otimes
 \Cal O((d+1)p))$ such that $\overline s=(s_1, \dots , s_{q_d})$ is a frame of
 $H^0(C,\overline{\varphi (g)}\otimes \Cal O(dp))$.    We have a
 commutative diagram:
 $$\matrix\format\c&&\quad\c\\
 0  &\longrightarrow  &H^0(C,\overline{\varphi (g)}\otimes \Cal O(dp))
 &\longrightarrow  &H^0(C,\overline{\varphi (g)}\otimes \Cal O((d+1)p))
 \\ \vspace{1\jot}
 && @VVV  @VVV \\ \vspace{1\jot}
 0  &\longrightarrow
 &H^0(C,{\widetilde{\Cal W}^o}_{d\vert C\times \varphi_d(s)})
 &\longrightarrow  &H^0(C,{\widetilde{\Cal W}^o}_{d\vert C\times \varphi_d(s)}
\otimes
 \Cal O(p)) ~, \endmatrix
 $$
 where the vertical maps are isomorphisms. Observe that, under the first
vertical isomorphism, the frame $\overline s$ is mapped to the frame
$\delta(\overline s)$. Now define
 $\alpha(s)$  to be the frame in $H^0(C,{\widetilde{\Cal W}^o}_{d\vert C\times
\varphi _d(s)} \otimes
 \Cal O(p))$           which is the image of the frame $s$ under  the second
vertical isomorphism.
Then $\alpha$ is a $P$-equivariant morphism.

 We claim that (as line bundles on $\Cal F_{d,Y}$)
 $$
 \varphi ^\ast_d
 (\widehat{\theta^\ast_d\frak L})
 \approx (\alpha^\ast \gamma^\ast_{d}
 (\widehat{\theta^\ast_{d+1}\frak L}))^U.
 \tag 3$$
 This follows since
 $$\align
 (\alpha^\ast &\gamma^\ast_{d}
 (\widehat{\theta^\ast_{d+1}\frak L}))^U \\
 &\approx \delta^\ast ((\gamma^\ast_{d}
 (\widehat{\theta^\ast_{d+1}\frak L}))^U) \\
 &\approx\varphi^\ast_d\sigma^\ast((\gamma^\ast_{d}
 (\widehat{\theta^\ast_{d+1}\frak L}))^U)\\
 &\approx \varphi^\ast_d(\widehat{\theta^\ast_{d}\frak L}), ~
 \text{using (2).}
 \endalign$$

 Now the bundle $\varphi^\ast_d
 (\widehat{\theta^\ast_{d}\frak L})$ has a GL$(q_d)$-linearization coming from
the action of $GL(q_d)$ on
 $\widehat{\theta^\ast_{d}\frak L}$ and (by definition of $\frak L_d$ ) the
bundle on $Y$
 obtained by descent is ${\frak L_d}_{\vert Y}$. On the
 other hand, the bundle $(\alpha^\ast\gamma^\ast_{d}
 (\widehat{\theta^\ast_{d+1}\frak L}))^U$
 has a GL$( q_d)$-action given by the action of $P/U \approx
 \text{GL}(q_d)$ arising from the action of $P$ on
 $\alpha^\ast\gamma^\ast_{d}
 (\widehat{\theta^\ast_{d+1}\frak L})$ (which in turn comes
 from the action of GL$(q_{d+1})$, in particular, $P$ on
 $\widehat{\theta^\ast_{d+1}\frak L})$ and the bundle on
 $Y$ obtained by descent via $\pi_d$ is ${\frak L_{d+1}}_{\vert Y}$.  Now
 on $\Cal F_{d,Y}^s := \pi_d^{-1}(Y^s)$, these two actions  of GL$( q_d)$
 coincide (i.e. the isomorphism  $\eta$ of line bundles on $\Cal F_{d,Y}$ as
guarnteed by (3) is GL$( q_d)$-equivariant on $\Cal F_{d,Y}^s$ ), as is seen
from the commutative diagram (got from the diagrams $D_1$ and $D_4$ ) (where
$Q^s_{d,Y} := \beta^{-1}_d(\Cal F^s_{d,Y})$ ):
 $$\matrix\format\c&&\quad\c \\
 &&Q^s_{d,Y}\\ \vspace{1\jot}
 &\sideset\beta_d\and\to\swarrow  && \sideset\and\gamma_{d}\circ
\alpha\to\searrow\\ \vspace{1\jot}
 \Cal F^s_{d,Y}  & \overset{\varphi_d}\to\longrightarrow  &R^s_d &R^s_{d+1}\\
\vspace{1\jot}
\pi_d  @VVV  @V \theta_d VV  \sideset\and\theta_{d+1}\to\swarrow\\
\vspace{1\jot}
 Y^s & \overset{\psi}\to\longrightarrow  &\frak M
 \endmatrix \tag $D_5$
 $$
 Since $Y^s$ is dense in $Y$, we have that $\Cal F^s_{d,Y}$ is dense in $\Cal
F_{d,Y}$;
 in particular, the isomorphism $\eta$ is GL$(q_d)$-equivariant on the whole of
$Y$.  Hence
 $
 {\frak L_d}_{\vert Y} \approx {\frak L_{d+1}}_{\vert Y}$. Denote this
isomorphism by $\mu$. Then the restriction of $\mu$ to $Y^s$ is the identity
map under the identification (1). From this it is easy to see that $
 {\frak L_d} \approx {\frak L_{d+1}}_{\vert X_d}$.
 This
 completes the proof  of the lemma.\qed
 \enddemo

 Finally we come to the
 \demo{{\bf (3.12)} Proof of Proposition (2.3)} For any Schubert variety
$X_\frak w$ , there exists a large enough $d(\frak w)$ such that
$X_\frak w \subset X_{d(\frak w)} .$  Let $\widehat{\frak L_\frak w}$
be the line bundle on $X_\frak w$ defined by  $\widehat{\frak L_\frak w} =
{\frak L_{d(\frak w)}}_{\vert X_{\frak w}} .$ By Lemma (3.11),  $\widehat{\frak
L_\frak w}$ is well defined and
$\widehat{\frak L_\frak w}_{\vert X^s_{\frak w}} \approx \psi^\ast (\frak
L)_{\vert X^s_{\frak w}}$,   where $\psi: X^s\to \frak M$ is the morphism as in
\S1.
  Moreover, for $\frak v \leq \frak w ~,  \widehat{\frak L_\frak w}_{\vert
X_{\frak v}} \approx  \widehat{\frak L_\frak v}~.$  In particular, by Lemma
(2.2), we get a line bundle  $\widehat{\frak L}$ on $X$ with  $\widehat{\frak
L}_{\vert X^s} \approx \psi^\ast (\frak L)$. This proves the proposition.
 \qed
 \enddemo
\vskip4ex

 \head 4.  Determination of  Pic $(\frak M)$
 \endhead

 \flushpar {\bf (4.1) Definition [D,\S2].}  Let $\frak{g}_1$ and
 $\frak{g}_2$ be two (finite dimensional) complex simple
 Lie algebras and $\varphi : \frak{g}_1 \to \frak{g}_2$
 be a Lie algebra  homomorphism.  There exists a unique
 number $m_\varphi \in \Bbb C$, called the {\it Dynkin index}
 of the homomorphism $\varphi$, satisfying
 $$
 \langle \varphi (x),\varphi (y) \rangle
 = m_\varphi \langle x,y \rangle, \text{ for all }x,y \in
 \frak{g}_1,
 $$
 where $\langle , \rangle$ is the Killing form on
 $\frak{g}_1$ (and $\frak{g}_2$) normalized so that
 $\langle \theta,\theta \rangle=2$ for the highest root
 $\theta$.

It is easy to see from \cite{KNR, Lemma 5.2} that for a
finite dimensional representation $V$ of $\frak g_1$ given by a
Lie algebra homomorphism $\varphi : \frak{g}_1 \to sl(V)$, we have
$m_\varphi = m_V$, where $m_V$ is as in \cite{KNR, \S5.1} and  $sl(V)$
is the Lie algebra of trace $0$ endomorphisms of $V$.

\vskip2ex

By taking a representation $V$ of $G_2$ such that $m_V \neq 0$, and
using \cite{KNR, Corollary 5.6}, the following proposition
follows easily.

 \proclaim{(4.2) Proposition}  Let $G_1,G_2$ be two
 connected  complex simple algebraic groups.  Then for any algebraic
 group homomorphism $\varphi: G_1\to G_2$, the induced
 map at the third homotopy group level
  $$\varphi _\ast : \pi_3 (G_1) \approx \Bbb Z
 \longrightarrow
 \pi_3 (G_2) \approx \Bbb Z$$
 is given by the multiplication via the Dynkin index
 $m_{d\varphi} $ of the induced Lie algebra  homomorphism
 $d\varphi :\frak{g}_1\to \frak{g}_2$, where $\frak{g}_1$ (resp. $\frak{g}_2$)
is the Lie
algebra of $G_1$ (resp.  $G_2$).

 In particular, $m_\varphi $ is an integer.
 \endproclaim
 \flushpar {\bf (4.3) Remark.}  The integrality of
 $m_{\varphi }$ is proved by Dynkin \cite{D, Theorem 2.2}, and so is the
 following lemma  \cite{D, Theorem 2.5}, by a quite different
 (and long) argument.

 \proclaim{(4.4) Lemma}  Let $\frak{g}$ be a complex simple Lie algebra and let
$V(\lambda)$ be  an
 irreducible representation of $\frak{g}$ with highest
 weight $\lambda$.  Then the Dynkin index
 $m_{V(\lambda)}$ of the representation $V(\lambda)$  is given by
 $$
 m_{V(\lambda)} =
 (\Vert \lambda+\rho \Vert^2 - \Vert \rho \Vert^2)
 \frac{\dim_{\Bbb C}V(\lambda) }{\dim_{\Bbb C}\frak{g}},
 $$
 where $\rho$ is the half sum of positive roots and the
 Killing form on $\frak{g}$ is normalized (as earlier) so
 that $\Vert \theta \Vert^2=2$ for the highest root
 $\theta$.
 \endproclaim

 \demo{Proof} The representation $V=V(\lambda)$ of course
 gives rise to a Lie algebra homomorphism $\varphi=\varphi _V:
 \frak{g}\to sl (V)$.  Since $m_V=m_{\varphi }$ (cf. \S4.1),  for any
 $x,y\in \frak{g}$
 $$
 m_V\langle x,y \rangle = \text{ trace } (\varphi
 (x)\circ \varphi (Y)).
 \tag1$$
 Choose a basis $\{ e_i\}$ of $\frak{g}$ and let $\{
 e^i\}$ be the dual basis of $\frak{g}$ with respect to
 the Killing form $\langle , \rangle$. Consider the
 Casimir element $\Omega:= \sum _i e_ie^i\in
 U(\frak{g})$.  Then $\Omega$ acts on $V$ via
 $$
 \Omega_V := \sum_i\varphi (e_i)\circ \varphi (e^i).
 \tag2$$
 But $V$ being irreducible of highest weight $\lambda$,
 $$
 \Omega_V = (\Vert \lambda+\rho \Vert^2 - \Vert \rho
 \Vert^2)
 I_V,
 \tag3$$
where $I_V$ is the identity operator of $V$.
 In particular,
 $$\align
 m_V  &= \frac{1}{\dim \frak{g}} \sum_i \text{
 trace }
 (\varphi (e_i) \circ \varphi (e^i)), \qquad \text{by
 (1)} \\
 &= \frac{1}{\dim \frak{g}} \text{ trace
 }\Omega_V,\qquad \text{by (2)} \\
 &= \frac{1}{\dim \frak{g}}
 (\Vert \lambda+\rho \Vert^2 - \Vert \rho \Vert^2) \dim V,
 \qquad \text{by (3)}.
 \endalign$$
 This proves the lemma. \qed
 \enddemo
\vskip1ex

 We also need the following
 \proclaim{(4.5) Lemma}  Let $\frak{g}$ be a complex
 simple Lie algebra and let $V$ and $W$ be two finite
 dimensional representations of $\frak{g}$.  Then
 $$
 m_{V\otimes W} = m_V\dim W + m_W \dim V.
 $$
 \endproclaim

 \demo{Proof} Write the characters
 $$\align
 \text{ch }V  &= \sum_\lambda n_\lambda e^\lambda  ~, \qquad
 \text{and} \\
 \text{ch }W  &= \sum_\mu m_\mu e^\mu~, \qquad \text{for
 some } n_\lambda,\ m_\mu \in \Bbb Z.
 \endalign$$
 Then
 $$
 \text{ch }(V\otimes W) = \sum_{\lambda,\mu }
 n_\lambda m_\mu e^{\lambda+\mu }.
 $$
 Hence by \cite{KNR, Lemma 5.2},
 $$\align
 2m_{V\otimes W} &=
 \sum_{\lambda,\mu } n_\lambda m_\mu
 \langle \lambda+\mu ,\theta^\vee \rangle^2 \\
 &= \sum n_\lambda m_\mu
 \langle \lambda,\theta^\vee \rangle^2
 + \sum n_\lambda m_\mu  \langle \mu,\theta^\vee \rangle^2+
 2  \sum n_\lambda m_\mu
 \langle \lambda,\theta^\vee \rangle  \langle \mu
 ,\theta^\vee\rangle\\
 &= 2\left( \sum_{\mu }m_\mu \right) m_V + 2\left(
 \sum_{\lambda }n_\lambda \right) m_W +
 2\left( \sum_\lambda n_\lambda \langle \lambda,\theta^\vee
  \rangle \right)
 \left( \sum_\mu m_\mu \langle \mu,\theta^\vee
 \rangle \right)\\
 &= 2(\dim W) m_V + 2(\dim V) m_W +
 2 \left(
 \sum_\lambda n_\lambda \langle \lambda,\theta^\vee
 \rangle \right)
 \left( \sum_\mu m_\mu \langle \mu,\theta^\vee
 \rangle \right) . \tag1
 \endalign$$

 For any $h\in \frak h$, define $\beta_V(h)= \sum_\lambda
 n_\lambda\langle \lambda,h \rangle$.  Then the map
 $\beta_V:\frak h\to \Bbb C,\ h\mapsto \beta_V(h)$ is
 $W$-equivariant (with the trivial action of $W$ on $\Bbb
 C$).  But $\frak{h}$ being an irreducible $W$-module,
 $$
 \beta_V \equiv 0.
 \tag2$$

 Combining (1) and (2), the lemma follows. \qed
 \enddemo

 \flushpar {\bf (4.6) Definition.}  Let $\frak{g}$ be a
 complex simple Lie algebra  and let $\theta$ be the
 highest root (with respect to some choice of the set of
 positive roots).  Express the associated coroot $\theta
 ^\vee$ in terms of the simple coroots:
 $$
 \theta^\vee = \sum^\ell _{i=1} m_i \alpha_i^\vee
 .
 $$

 Now define $d=d(\frak{g})$ to be the least common
 multiple of $\{ m_i\}_{i=1,\dots ,\ell }$.  Then the
 number $d$ is given as follows:

 $$
 \matrix\format\c&\quad\c\qquad\qquad\qquad &\c\\
 \text{{\bf Type of }}\frak g &&\text{{ d}}\bold( \frak g\bold) \\ \\
 A_\ell (\ell \geq 1),  &C_\ell (\ell \geq 2) &1\\
 B_\ell  &\ \ (\ell \geq 3) &2\\
 D_\ell  &\ \ (\ell \geq 4) &2\\
 G_2 &&2\\ F_4&&6\\ E_6&&6\\ E_7&&12\\E_8&&60
 \endmatrix
 $$

 \proclaim{(4.7) Proposition}  For any finite dimensional
 representation $V$ of $\frak{g}$, the number
 $d(\frak{g})$ divides $m_V$.  Moreover, there exists an
 irreducible representation $V_o$ of $\frak{g}$ such that
 $d(\frak{g})= m_{V_o}$.
 \endproclaim

 \demo{Proof} Unfortunately, our proof is case by case.
 We follow the indexing convention as in \cite{B, Planche
 I-IX}.  We denote the $i$-th fundamental weight $( 1\leq i\leq
 \ell )$ by $\omega_i$.
\vskip1ex

 \flushpar {\bf Case 1} -- $A_\ell (\ell \geq 1),\ C_\ell
 (\ell \geq 2)$:~
 As in \cite{KNR, Lemma 5.2}, $m_{V_o}=1$, for the
 standard $(\ell +1)$-dimensional representation $V_o$ of
 $A_\ell $.  Similarly for the standard $2\ell$-dimensional
 representation $V_o$ of $C_\ell $ (with
 highest weight $\omega_1$), $m_{V_o}=1$ (as can be seen
 from Lemma 4.4).
 \vskip1ex

   For a simply-connected group $G$, since the
 fundamental representations $\{ V(\omega_i) \}_{1\leq
 i\leq \ell }$ generate the representation ring $R(G)$
 as an algebra (cf \cite{A}), to prove that $d(\frak{g})$
 divides $m_V$ for any $\frak{g}$-module $V$, it suffices
 to show that $d(\frak{g})$ divides $m_i:=m_{V(\omega_i)}$
 for all $1\leq i\leq \ell $ (cf. Lemma 4.5).  In the
 following calculations, we make use of Lemma (4.4)
 and \cite{B, Planche I-IX} freely.
 \vskip1ex

 \flushpar {\bf Case 2} -- $B_\ell\  (\ell \geq 3)$:
 ~For $1\leq i\leq \ell -1$, $m_i=2\bigl( \matrix 2\ell
 -1\\i-1\endmatrix \bigr)$, since dim $V(\omega_i)=\bigl( \matrix
 2\ell+1\\i\endmatrix \bigr)$; and $m_\ell=2^{\ell-2}$.

 In particular, $m_1=2$, so take $V_o=V(\omega_1)$.
 \vskip1ex
 \flushpar {\bf Case 3} -- $D_\ell \ (\ell \geq 4)$ :~
 For $1\leq i\leq \ell -2$, $m_i=2 \bigl( \matrix 2\ell
 -2\\i-1\endmatrix \bigr)$, since dim $V(\omega_i)= \bigl( \matrix
 2\ell \\i\endmatrix \bigr)$; and
 $m_{\ell-1} = m_\ell =2^{\ell
 -3}$.

 In particular, $m_1=2$.
 \vskip1ex

 In the following  calculations, dim$\,V(\omega_i)$ is
 taken from \cite{BMP}.

 \flushpar {\bf Case 4} -- $G_2$:  $m_1=2,\ m_2=8$.

 (Observe that $V(\omega_2)$ is the adjoint
 representation of $G_2$ and hence $m_2$ can be

  calculated from \cite{KNR, Lemma 5.2 and Remark 5.3}.)
 \vskip1ex
 \flushpar {\bf Case 5} -- $F_4$:  $m_1,\dots ,m_4$ are
 respectively 18, 9$\times$98, l26, and 6.

 \vskip1ex
 \flushpar {\bf Case 6} -- $E_6$:  $m_1,\dots ,m_6$ are
 respectively 6, 24, 150, 1800, 150, and 6.

 \vskip1ex
 \flushpar {\bf Case 7} -- $E_7$:  $m_1,\dots ,m_7$ are
 respectively 36, 360, 65$\times$72, 2750$\times$108,
 104$\times$165,

 8$\times$81, and 12.

 \vskip1ex
 \flushpar {\bf Case 8} -- $E_8$:  $m_1,\dots ,m_8$ are
 respectively 12$\times$125, 4750$\times$18,
 49$\times$108000, 75$\times$111275472,

 30$\times$4720170, 45$\times$39520, 15$\times$980, and
 60. \qed
 \enddemo
 \vskip2ex

 \flushpar {\bf (4.8) Remark.}  The values of $m_i$ given
 above are also contained in \cite{D}, but many of his
 values are incorrect.

 \vskip1ex

 Combining Proposition (4.7) and Theorem (2.4)
 with the chart in Definition (4.6), we get the following  strengthening of
Theorem
 (2.4).
 \proclaim{(4.9) Theorem}  With the notation and
 assumptions as in Theorem {\rm (2.4)}, consider the
 injective map $\overline{\psi^\ast}:\text{{\rm Pic}}\,(\frak
 M)\hookrightarrow \text{{\rm Pic}}\, (X) \approx \Bbb
 Z.$  Then
\vskip1ex
 \flushpar
 (1) $\overline{\psi^\ast}$ is surjective in the case
 where $G$ is of type $A_\ell ~ (\ell \geq 1),$   and $C_\ell ~ (\ell \geq 2).$

 \flushpar
(2) The order $\gamma =\gamma _G$ of the cokernel of
 $\overline{\psi^\ast}$ is bounded as follows :
\vskip1ex

{\rm (a)} $G=B_\ell \
 (\ell \geq 3), \quad \gamma \leq 2$

{\rm (b)}
 $ G= D_\ell \ (\ell \geq 4), \quad \gamma \leq 2$

{\rm (c)} $G=G_2, \quad \gamma \leq 2$

 {\rm (d)} $\ G=F_4, \quad \gamma \leq 6$

 {\rm (e)} $\ G=E_6, \quad \gamma \leq 6$

 {\rm (f)} $\ G=E_7, \quad \gamma \leq 12$

 {\rm (g)} $\ G=E_8, \quad \gamma \leq 60 ~.$
 \endproclaim

\flushpar {\bf (4.10) Definition.}  A (complex) representation $V$ of $G$ is
said to be {\it orthogonal} if there exists a $G$-invariant non-degenerate
symmetric ${\Bbb C}$-bilinear form on $V$.

Clearly, an orthogonal representation is isomorphic with its dual. Of course
the adjoint representation ${\frak g}$ of $G$ is orthogonal.

Even though we do not have a full proof as yet (nor do we know any place in the
literature where it is proved), we believe that the following proposition is
true. (G. Faltings has written to the first author that this should be true. He
has suggested that  one should  show that the cohomology is locally
representable by a complex $\Cal E \overset S \to \longrightarrow \Cal E^\ast
$, where $S$ corresponds to an alternating form.  Thus giving det $(S) \approx
$ Pfaff $(S)^2$.

 \proclaim{(4.11) ``Proposition"}   Let $V$ be an orthogonal representation of
$G$. Then the theta bundle $\Theta (V)$ on ${\frak M}$ admits a square root as
a line bundle,  i.e., there exists an
 (algebraic) line bundle $\frak L$ on $\frak M$ such that
 $\frak L^2 \approx \Theta (V)$.
 \endproclaim
\vskip2ex

  \flushpar {\bf (4.12) Remark.} The validity of the above proposition will
show that the map $\overline{\psi^\ast}$ (as in the above theorem 4.9) is
surjective for $G=B_\ell \
 (\ell \geq 3),
  G= D_\ell \ (\ell \geq 4),$ and $G=G_2$.
(Observe that for $G_2$ , $V(\omega_1)$ is seven dimensional orthogonal
representation.) Moreover
the bounds for $\gamma$ in the cases $F_4, E_7,$ and
 $E_8$ can be improved to $3, 6,$ and $30$ respectively.

\vskip2ex

We feel that the following question has an affirmative answer.

\flushpar {\bf (4.13) Question.}  For any $C,G$ as in
 Theorem (4.9) (with genus $C\geq 2$), is it true that
 the bounds for the  order of
 $\gamma $ as in Remark (4.12)  are
 in fact equalities?
 \vskip1ex

 For $C=\Bbb P^1(\Bbb C)$, since the moduli space $\frak
 M$ reduces to a point, the map $\overline{\psi^\ast}$ is
 clearly an isomorphism.


\vskip4ex
 \Refs \widestnumber\key{BMP}

 \ref
 \key A
 \by Adams, J. F.
 \book `Lectures On Lie Groups'
 \publ W. A. Benjamin, Inc.
 \yr 1969
 \endref



 \ref
 \key B
 \by Bourbaki, N.
 \book `Groupes et Alg\`ebres de Lie
 \bookinfo Chap. IV--VI'
 \publ Hermann
 \yr 1968
 \publaddr Paris
 \endref


 \ref
 \key Be
 \by Behrend, K. A.
 \paper Semi-stability of reductive group schemes over
 curves
 \jour Math. Annalen
 \vol 301
 \yr 1995
 \pages 281--305
 \endref


 \ref
 \key BMP
 \by Bremner, M. R.; Moody, R. V. and Patera, J.
 \book `Tables of Dominant Weight Multiplicities for
 Representations of Simple Lie Algebras'
 \yr
 \endref

\ref
\key Bo
\by Boutot, J. F.
\paper Singularites rationelles et quotients par les groupes reductifs
\jour Invent. Math.
\vol 88
\yr 1987
\pages 65-68
\endref

 \ref
 \key D
 \by Dynkin, E. B.
 \paper Semisimple subalgebras of semisimple Lie algebras
 \jour Amer. Math. Soc. Transl., Ser. II,
 \vol 6
 \yr 1957
 \pages 111--244
 \endref

\ref
 \key DN
 \by Drezet, J.-M. and  Narasimhan, M. S.
 \paper Groupe de Picard des vari\'et\'es de modules de
 fibr\'es
 semi-stables sur les courbes alg\'ebriques
 \jour Invent. Math.
 \vol 97
 \yr 1989
 \pages 53--94
 \endref

\ref
\key F
\by Faltings, G.
\paper Stable $G$-bundles and projective connections
\jour J. Alg. Geom.
\vol 2
\yr 1993
\pages 507-568
\endref

 \ref
 \key GR
 \by Grauert, H. and Riemenschneider, O.
 \paper Verschwindungss\"atze f\"ur analytische
 cohomologiegruppen auf komplexen Ra\"umen
 \inbook Several Complex Variables I
 \publ Springer Lecture Notes
 \publaddr Maryland
 \yr 1970
 \endref

 \ref
 \key H
 \by Hartshorne, R.
 \book `Algebraic Geometry'
 \publ Springer Verlag
 \publaddr Berlin-Heidelberg-New York
 \yr 1977
 \endref

 \ref
 \key K
 \by Kac, V. G.
 \book `Infinite Dimensional Lie Algebras'
 \bookinfo Third Edition
 \publ Cambridge University Press
 \yr 1990
 \endref

\ref
\key Ku
\by Kumar, S.
\paper Symmetric and exterior powers of homogeneous vector bundles
\jour Math. Annalen
\vol 299
\yr 1994
\pages 293--298
\endref

\ref
 \key KNR
 \by Kumar, S.; Narasimhan, M. S. and Ramanathan, A.
 \paper Infinite Grassmannians and moduli spaces of
 $G$-bundles
 \jour Math. Annalen
 \vol 300
 \yr 1994
 \pages 41--75
 \endref

 \ref
 \key Kr
 \by Kraft, H.
 \paper Algebraic automorphisms of affine space
 \inbook Proceedings of the Hyderabad Conference on
 Algebraic Groups
 \bookinfo ed. S. Ramanan
 \publ Manoj Prakashan
 \yr 1991
 \pages 251--274
 \endref

 \ref
 \key NR
 \by Narasimhan, M. S. and Ramadas, T. R.
 \paper Factorisation of generalised theta functions. I
 \jour Invent. Math.
 \vol 114
 \yr 1993
 \pages 565--623
 \endref


\ref
\key PV
\by    Verdier, J.L.  and  Potier, J. Le
\paper
 Modules des fibr\`es stables sur les courbes
 alg\`ebriques, expos\'e 4
\jour   Progress in Mathematics
\vol
 54
\publ Birkhauser
\yr 1985
\endref

\ref
\key R1
\by Ramanathan, A.
\book `Stable principal bundles on  a compact Riemann surface- Construction of
moduli space'
\bookinfo Thesis (Bombay University, Bombay, India)
\yr 1976
\endref

\ref
\key R2
\by Ramanathan, A.
\paper Moduli of principal bundles
\jour  L. N.
 M.
\vol 732
\publ ( Algebraic Geometry, Proceedings Copenhagen)
 \yr 1978
\pages 527--533
\endref

\endRefs
\enddocument